\documentclass[preprint]{elsarticle}
\usepackage{amssymb}
\usepackage{graphicx}
\usepackage{geometry}
 \usepackage{subfigure}
  \usepackage{subfigmat}
\usepackage{indentfirst}
\usepackage{amsmath}
\usepackage{textcomp}
\usepackage{amsthm}
\usepackage{mathrsfs}
\usepackage{hyperref}
\usepackage{flafter}
\usepackage{sidecap}
\usepackage{float}

\newtheorem{proposition}{Proposition}

\author[rvt]{Hantian Zhang\corref{cor1}\fnref{fn1}}
\ead{h.zhang.ae.mech@hit.edu.cn}
\address[rvt]{School of Astronautics, Harbin Institute of Technology, Harbin, China, 150001}
\author[focal]{ Dong Eui Chang\fnref{fn2}}
\ead{dechang@math.uwaterloo.ca}
\address[focal]{Department of Applied Mathematics, University of Waterloo, Waterloo, ON, Canada, N2L 3G1}
\author[rvt]{Qingjie Cao\fnref{fn3}}
\ead{q.j.cao@hit.edu.cn}
\fntext[fn1]{Undergraduate Student; h.zhang.ae.mech@hit.edu.cn} 
\fntext[fn2]{Associate Professor; dechang@math.uwaterloo.ca}
\fntext[fn3]{Professor; q.j.cao@hit.edu.cn}

\title{Lyapunov-based Low-thrust Optimal Orbit Transfer:\\ An approach in Cartesian coordinates}

\cortext[cor1]{Corresponding author}

\begin{document}

\begin{abstract}
This paper presents a simple approach to low-thrust optimal-fuel and optimal-time transfer problems between two elliptic orbits using the Cartesian coordinates system. In this case, an orbit is described by its specific angular momentum and Laplace vectors with a free injection point. Trajectory optimization with the pseudospectral method and nonlinear programming are supported by the initial guess generated from the Chang-Chichka-Marsden Lyapunov-based transfer controller. This approach successfully solves several low-thrust optimal problems. Numerical results show that the Lyapunov-based initial guess overcomes the difficulty in optimization caused by the strong oscillation of variables in the Cartesian coordinates system. Furthermore, a comparison of the results shows that obtaining the optimal transfer solution through the polynomial approximation by utilizing Cartesian coordinates is easier than using orbital elements, which normally produce strongly nonlinear equations of motion. In this paper, the Earth's oblateness and shadow effect are not taken into account.
\end{abstract}

\maketitle

\section{Introduction}

Three dimensional low-thrust optimal orbit transfer has attracted much inquiry focused on trajectory optimization using orbital elements. Due to  strong nonlinearity of differential equations in Gaussian form with orbital elements, it is often difficult to obtain the optimal solution numerically in this system.

Some of the earliest works on the orbit transfer between neighboring elliptic orbits and on the transfer between coplanar and coaxial ellipses were presented by Edelbaum \cite{edelbaum1964optimum, edelbaum1965optimum}. However, his elements as well as the Keplerian  elements all contain singularity. To avoid the singularity, Arsenault \cite{arsenault1970} firstly introduced the equinoctial elements. Broucke and Cefola \cite{broucke1972equinoctial} developed nonsingular equinoctial orbital elements using the Lagrange and Poisson brackets of Keplerian elements. Kechichian \cite{kechichian1996optimal} presented an application of these nonsingular elements to solve the minimum time rendezvous problem with constant acceleration. Chobotov \cite{chobotov2002orbital} considered more cases in minimum time transfer, including the comparison between exact solutions and approximate solutions obtained by the averaging technique. Gerffroy and Epenoy \cite{geffroy1997optimal} made further progress in both minimum time and fuel transfer problems using the averaging technique with the constraints of Earth's oblateness and shadow effect taken into account. More recent works using the numerical averaging technique were presented by Tarzi and Speyer et al.\cite{tarzi2013fuel}, who provided a quick and accurate numerical approach for a wide range of transfers, including orbital perturbations such as Earth's oblateness and shadow effect. Besides the strong nonlinearity, an additional disadvantage in using the equinoctial elements is the complexity of the equations in this coordinate system. With this system, Kepler's equation must be solved by iteration to get the eccentric longitude at each integration step. Hence, equinoctial elements present challenges for trajectory optimization. Walker \cite{walker1985set,walker1986set} put forth another important development in the study of equinoctial elements when he used the Stroboscopic method to modify orbital elements. He also altered differential equations into an approximative form, containing five dependent variables and one independent variable, as a means of achieving faster computation performance without solving Kepler's equation. Roth \cite{roth1979higher} introduced the Stroboscopic method to obtain a higher order approximation for small perturbed dynamical systems, which depends on several slow variables and one fast variable. Haberkorn and Gergaud \cite{haberkorn2004low,gergaud2006homotopy} investigated the application of the modified orbital elements by using the homotopy method. Cui et al. \cite{cui2008low} used sequential quadratic programming under their Lyapunov feedback control law, which was based on a function made up of modified elements, to obtain the optimal-Lyapunov solution without optimal transfer. Gao and Li \cite{gao2010optimization} made the similar work to optimize the Lyapunov function but never reached the optimal solution based on their Lyapunov control law.

An advanced technique using the Lyapunov-based controller to solve the low-thrust orbit transfer problem in Cartesian coordinates was introduced and rigorously proved by Chang et al. \cite{chang2002lyapunov}. This technique is based on the fact that a non-degenerate Keplerian orbit is uniquely described by its specific angular momentum and Laplace vectors. The resulting Lyapunov function provides an asymptotically stabilizing feedback controller, such that the target elliptic Keplerian orbit becomes a locally asymptotically stable periodic orbit. However, the Lyapunov-based transfer trajectory is not optimal in every sense. In this paper, the Lyapunov-based transfer presented in \cite{chang2002lyapunov} shall be called Chang-Chichka-Marsden (CCM)\footnote{Abbreviation for Chang-Chichka-Marsden} transfer to distinguish it from any other Lyapunov-based transfers. 

The motivation behind this paper is to use the CCM transfer trajectory as an initial guess for optimization in order to obtain the optimal transfer solution utilizing Cartesian coordinates. Using this method avoids the numerical disadvantages due to strong nonlinearity and complexity in the use of orbital elements.

This paper reviews the CCM transfer method \cite{chang2002lyapunov}  in Section 2. In Section 3, a means to translate Keplerian elements into specific angular momentum and Laplace vectors is presented. Section 4 presents the proposed approach and the optimality (KKT)\footnote{Abbreviation for Karush-Kuhn-Tucker} conditions for the minimum fuel consumption problem in Cartesian coordinates. Specifically, the Chebyshev-Gauss pseudospectral method is introduced to illustrate how the continuous optimal control problem can be reduced to a discretized nonlinear programming problem. Finally, in Section 5, numerical simulations are carried out to make detailed comparisons between the optimal results using Cartesian coordinates and those using orbital elements with the same initial guess. It shows that the use of Cartesian coordinates makes it easier to obtain the correct optimal solution.

\section{Chang-Chichka-Marsden Transfer}
This section summarizes and reviews the Chang-Chichka-Marsden (CCM) transfer in \cite{chang2002lyapunov}. This transfer employs Lyapunov-based controllers to achieve asymptotically stable transfers between elliptic orbits in a two-body problem. 
\subsection{Two-Body Problem}
This paper assumes that the configuration space is $\mathbb R^{3}_0:=\mathbb R^3-\{(0,0,0)\}$. Let $T\mathbb R^{3}_0=(\mathbb R^{3}-\{(0,0,0)\})\times \mathbb R^3$ be the tangent space of $\mathbb R^{3}_0$,  and $(\mathbf r, \dot{\mathbf r})$ be the coordinates on $T\mathbb R_0^3$. The equations of motion are given by
\begin{equation}
\label{eq:twobody}
\ddot{\mathbf r}=-\mu\frac{\mathbf r}{{\| \mathbf r \|}^3}
\end{equation}
of which the solutions are regarded as the \textsl{Keplerian orbits}, where $\mu$ is the gravitational constant. The specific energy $\varepsilon: T\mathbb R^{3}_0 \to \mathbb R$ is defined by
\[
\varepsilon(\mathbf r, \dot{\mathbf r})=\frac{1}{2} {\| \dot{\mathbf r} \|}^2 - \frac{\mu}{\|{\mathbf{r}}\|}.
\]
Define $\pi=(\mathbf L, \mathbf A):  T\mathbb R^{3}_0 \to \mathbb R^3 \times \mathbb R^3$ by
\begin{align*}
\mathbf L(\mathbf r, \dot{\mathbf r})&=\mathbf r \times \dot{\mathbf r},\\
\mathbf A(\mathbf r, \dot{\mathbf r})&={\dot{\mathbf r}\times (\mathbf r \times \dot{\mathbf r})}-{\mu}\frac{\mathbf r}{\| \mathbf r \|},
\end{align*}
where $\mathbf L$ is the specific angular momentum vector and $\mathbf A$ is the Laplace vector. The Laplace vector $\mathbf A$ is related to the eccentricity vector $\mathbf e$ as follows:
\[
\mathbf A = \mu \mathbf e.
\]
The three quantities $\varepsilon$, $\mathbf L$ and $\mathbf A$ satisfy the following two identities:
\[
\mathbf L \cdot \mathbf A = 0, \quad \| \mathbf  A\|^2 = \mu^2 + 2\varepsilon \| \mathbf L\|^2. 
\]
Define the sets
\begin{align}
\Sigma_e &=\{ (\mathbf r, \dot{\mathbf r})\in T\mathbb R^{3}_0 \mid\mathbf L(\mathbf r, \dot{\mathbf r}) \neq 0, \| \mathbf A (\mathbf r, \dot {\mathbf r})\| < \mu  \}, \nonumber \\
D &=\{ (\mathbf x, \mathbf y)\in \mathbb R^3 \times \mathbb R^3 \mid \mathbf x \cdot \mathbf y =0, \mathbf x \neq 0, \| \mathbf y \| < \mu \}. \label{set:D}
\end{align}
The following Proposition is from \cite{chang2002lyapunov}.
\begin{proposition}
The following hold:\\
1. $\Sigma_e$ is the union of all non-degenerate elliptic Keplerian orbits.\\
2. $\pi(\Sigma_e) = D$ and $\Sigma_e = \pi^{-1}(D)$.\\
3. The fiber $\pi^{-1}(\mathbf x, \mathbf y)$ consists of a unique (oriented) non-degenerate elliptic Keplerian orbit for each $ (\mathbf x, \mathbf y) \in D$.
\label{prop1}
\end{proposition}

\subsection{Chang-Chichka-Marsden  Transfer}
The equation of the motion with a specific control force is given by
\begin{equation}
\label{eq:newton}
\ddot{\mathbf r}=-\mu\frac{\mathbf r}{{\| \mathbf r \|}^3}+\mathbf F,
\end{equation}
where $\mathbf F$ is the control force. 
Define a metric $d_k$ on $\mathbb R^3 \times \mathbb R^3$ by
\begin{displaymath}
d_k((\mathbf x_1,\mathbf y_1),(\mathbf x_2,\mathbf y_2))=\sqrt{\frac{1}{2}k{\|\mathbf x_1-\mathbf x_2\|}^2+\frac{1}{2}\|\mathbf{y_1-\mathbf y_2\|}^2}
\end{displaymath}
with $k > 0$ an arbitrary parameter, and $(\mathbf x_1,\mathbf y_1),(\mathbf x_2,\mathbf y_2) \in \mathbb R^3 \times \mathbb R^3$. Let $B_{d_k}((\mathbf x, \mathbf y), r) \subset \mathbb R^3 \times \mathbb R^3 $ be the open ball of radius $r$ centered at $(\mathbf x,\mathbf y) \in \mathbb R^3 \times \mathbb R^3$ in the $d_k$-metric and $\bar B_{d_k} ((\mathbf x, \mathbf y), r)$  its closure.

Let $(\mathbf L_T , \mathbf A_T ) \in D$ be the pair of the angular momentum and the eccentricity vector of a target elliptic orbit. Define a Lyapunov function $V$ on $T\mathbb R^{3}_{0}$ by
\begin{equation}
\label{eq:lyapunov}
V(\mathbf r, \dot{\mathbf r})=\frac{1}{2}k\| \mathbf L(\mathbf r, \dot{\mathbf r})-\mathbf L_T \|^2+\frac{1}{2}\| \mathbf A(\mathbf r, \dot{\mathbf r})-\mathbf A_T \|^2.
\end{equation}
Along the trajectory of (\ref{eq:newton}) there is
\begin{align*}
\frac{d}{dt}\mathbf L(\mathbf r, \dot{\mathbf r})&= \mathbf r \times \mathbf F, \\
\frac{d}{dt}\mathbf A(\mathbf r, \dot{\mathbf r})&= \mathbf F \times \mathbf L(\mathbf r, \dot{\mathbf r})+\dot{\mathbf r} \times (\mathbf r \times \mathbf F).
\end{align*}
Hence,
\begin{displaymath}
\frac{d}{dt}V(\mathbf r,\dot{\mathbf r})=\mathbf F \cdot \Big ( k\Delta\mathbf L \times \mathbf r + \mathbf L(\mathbf r, \dot{\mathbf r}) \times \Delta \mathbf A + (\Delta\mathbf A \times \dot{\mathbf r}) \times \mathbf r  \Big ),
\end{displaymath}
where
\begin{equation}
\Delta \mathbf L=\mathbf L(\mathbf r, \dot{\mathbf r}) - \mathbf L_T \quad, \quad \Delta \mathbf A = \mathbf A(\mathbf r, \dot{\mathbf r})-\mathbf A_T.
\end{equation}
Take a controller as follows:
\begin{equation}
\label{eq:control}
\mathbf F(\mathbf r, \dot{\mathbf r}; \mathbf L_T, \mathbf A_T)=-f(\mathbf r, \dot{\mathbf r})\left( k\Delta\mathbf L \times \mathbf r +  \mathbf L(\mathbf r, \dot{\mathbf r}) \times \Delta \mathbf A + (\Delta\mathbf A \times \dot{\mathbf r}) \times \mathbf r \right)
\end{equation}
with an arbitrary function $f(\mathbf r, \dot{\mathbf r})>0$. Then,
\begin{equation}
\frac{d}{dt}V(\mathbf r,\dot{\mathbf r})=-f(\mathbf r, \dot{\mathbf r})\left\| k\Delta\mathbf L \times \mathbf r + \mathbf L(\mathbf r, \dot{\mathbf r}) \times \Delta \mathbf A + (\Delta\mathbf A \times \dot{\mathbf r}) \times \mathbf r  \right\|^2 \leq 0.
\end{equation}
The following Proposition is proven  in \cite{chang2002lyapunov} using LaSalle's invariance principle \cite[pp.58-59]{la1961stability}.  

\begin{proposition}\label{prop:2}
Let $(\mathbf L_T,\mathbf A_T) \in D$ be the pair of the specific angular momentum and Laplace vectors of the target elliptic orbit. Take any closed ball $\bar B_{d_k}((\mathbf L_T, \mathbf A_T),l)$ of  radius $l>0$ centered at $(\mathbf L_T, \mathbf A_T)$ contained in the set $D$ defined in (\ref{set:D}).
Then every trajectory starting in the subset $\pi^{-1}(B_{d_k}((\mathbf L_T, \mathbf A_T),l))$ of $T\mathbb R^3_0$ remains in that subset and asymptotically converges to the target elliptic orbit $\pi^{-1}(\mathbf L_T, \mathbf A_T)$ in the closed-loop dynamics (\ref{eq:newton}) with the control law in (\ref{eq:control})
\end{proposition}

The choice of the parameter $k$ in the Lyapunov function $V$  plays a crucial role in determining the transfer trajectory \cite{chang2002lyapunov}. It determines the relative weighting between the two quadratic terms in the function $V$ in (\ref{eq:lyapunov}). With a small $k$  the shape of the trajectory will adjust to that of the target orbit first because a more weight is on $\|\mathbf A\|^2$. On the other hand, with  a large $k$, the normal direction  of the trajectory plane will adjust to that of the target orbit plane first because a more weight is on $\|\mathbf L\|^2$. The parameter $k$ also determines the shape of the region of attraction since $k$ determines the shape of the ball $B_{d_k}$ in the metric $d_k$. Additionally, the CCM transfer works well for parabolic transfer, although the success of the transfer is proven exclusively  for  elliptic orbits only in Proposition \ref{prop:2}.

\section{Transform of Orbital Elements}
This section presents the transform of the six Keplerian  elements to specific angular momentum and Laplace vectors for convenient reference. The state vector at periapsis to derive the transform in Cartesian coordinates with the Earth at the origin.  Let
\[
h = \| \mathbf L \|, \quad e=\|\mathbf e\|.
\]
The periapsis of the orbit in the geocentric equatorial frame is determined by
\begin{displaymath}
(\mathbf r, \dot{\mathbf r})=\left(\frac{h^2}{\mu (e+1)} \frac{\mathbf A}{\| \mathbf A \|},  \frac{\mu (e+1)}{h} \frac{\mathbf L}{\| \mathbf L \|} \times \frac{\mathbf A}{\| \mathbf A \|} \right).
\end{displaymath}
In the perifocal frame,
\begin{displaymath}
\{ \mathbf r \}_{P}=\frac{h^2}{\mu (e+1)} \left\{ \begin{array}{ll}1\\0\\0 \end{array} \right\} \quad , \quad
\{ \dot{\mathbf r} \}_{P}=\frac{\mu(e+1)}{h}\left\{ \begin{array}{ll} 0\\1\\0 \end{array} \right\}.
\end{displaymath}
The transformation matrix \cite[p.174]{curtis2005orbital} from the perifocal frame into the geocentric equatorial frame is given by 
\begin{displaymath}
[\mathbf Q]_{PE}=\left[ \begin{array}{ccc}
\cos\Omega \cos\omega - \sin\Omega \sin\omega\cos i & -\cos\Omega \sin\omega -\sin\Omega \cos i \cos\omega & \sin\Omega \sin i \\
\sin\Omega \cos\omega + \cos\Omega \cos i \sin\omega & -\sin\Omega \sin\omega+\cos\Omega \cos i \cos\omega & -\cos\Omega \sin i \\
\sin i \sin\omega & \sin i \cos\omega & \cos i
\end{array} \right].
\end{displaymath}
The state vector in the geocentric equatorial frame is found by carrying out the matrix
multiplications
\begin{displaymath}
\{\mathbf r\}_E=[\mathbf Q]_{PE}\{\mathbf r\}_P, \{\dot{\mathbf r}\}_E=[\mathbf Q]_{PE}\{\dot{\mathbf r}\}_P.
\end{displaymath}
Thus, the components of $\frac{\mathbf A}{\| \mathbf A \|}$ and $\frac{\mathbf L}{\| \mathbf L \|} \times \frac{\mathbf A}{\| \mathbf A \|}$ are derived. Then using the identity $\frac{\mathbf L}{\| \mathbf L \|}=\frac{\mathbf A}{\| \mathbf A \|} \times \left(\frac{\mathbf L}{\| \mathbf L \|} \times \frac{\mathbf A}{\| \mathbf A \|}\right)$, the specific angular momentum and Laplace vectors in the geocentric equatorial frame are computed as follows:
\begin{align*}
{\mathbf L}=\sqrt{\mu a (1- e^2)}\left\{ \begin{array}{ccc} \sin i \sin\Omega \\ -\sin i \cos\Omega \\ \cos i \end{array} \right\},
\end{align*}
\begin{align*}
{\mathbf A}=\mu e \left\{ \begin{array}{ccc}
\cos\Omega \cos\omega-\sin\Omega \sin\omega \cos i \\
\sin\Omega \cos\omega+\cos\Omega \sin\omega \cos i \\
\sin i \sin\omega \end{array} \right\}.
\end{align*}
On equatorial orbits $(\Omega=0, i=0)$, they simplify to 
\begin{align*}
\mathbf L =\sqrt{\mu a (1-e^2)}\left\{ \begin{array}{ccc} 0 \\ 0 \\ 1 \end{array} \right\} \quad , \quad \mathbf A =\mu e \left\{ \begin{array}{ccc} \cos\omega \\ \sin\omega \\ 0 \end{array} \right\}.
\end{align*}
On circular orbits $(\omega=0, e=0)$, they  simplify to
\begin{align*}
{\mathbf L}=\sqrt{\mu a}\left\{ \begin{array}{ccc} \sin i \sin\Omega \\ -\sin i \cos\Omega \\ \cos i \end{array} \right\} \quad , \quad \mathbf A=0\left\{ \begin{array}{ccc} \cos\Omega \\ \sin\Omega \\ 0 \end{array} \right\}.
\end{align*}

\section{Optimal Orbit Transfer}

The CCM transfer trajectory is used as an initial guess to support the trajectory optimization in the open-loop system using the direct Chebyshev-Gauss pseudospectral transcription method and a nonlinear programming solver. The Cartesian coordinates and the modified orbital elements in optimization are compared. 

\subsection{Optimization in Cartesian Coordinates}
Let $\mathbf x=(\mathbf r, \dot{\mathbf r})$ denote the Cartesian coordinates in the geocentric equatorial frame. The minimum fuel consumption problem is given as follows:
\begin{displaymath}
(P1)  \left \{ \begin{array}{lll} \text{Minimize} &  J=\int^{t_f}_{t_0} \| \mathbf u(t) \|dt \\ 
\text{Subject to} & \dot{\mathbf x}(t)= \mathbf m(t) \mathbf x(t) + \mathbf u(t)\,  , \,  \forall t \in [t_0,t_f] \\
 & \| \mathbf u(t) \| \leq F_{\rm max} \, , \, \forall t \in [t_0,t_f] \\
& \mathbf x(t_0) \, \text{ fixed}\\
& \mathbf L_T, \mathbf A_T \, \text{ fixed}\\
\end{array}  \right.
\end{displaymath}
with
\begin{displaymath}
\mathbf m=\left[  \begin{array}{cc}
0 & \mathbf E_3\\
d\mathbf E_3 & 0\\
\end{array}  \right]_{6 \times 6}, \quad  \quad 
d=-\frac{\mu}{\| \mathbf r \|^3}, \quad  \quad \mathbf u=[0,0,0,F_1,F_2,F_3]^{T},
\end{displaymath}
where  $\mathbf x$ denotes the state vector, $\mathbf E_3$ is the $3\times 3$ identity matrix and $\mathbf u$ the control vector. The boundary conditions are given by
\begin{align}
& \mathbf S_0(\mathbf x(t_0))=\mathbf x(t_0) - \mathbf x_0=\mathbf 0 \\
& \mathbf S_{L}(\mathbf x(t_f))=\mathbf L(\mathbf x(t_f))-  \mathbf L_{T}= \mathbf 0 \\
& \mathbf S_{A}(\mathbf x(t_f))= \mathbf A(\mathbf x(t_f))-  \mathbf A_{T}= \mathbf 0
\end{align}
For the minimum time problem, the cost function  $J=\int_0^{t_f}dt$ shall be used.


To reduce the continuous optimal control problem (OCP) into a discretized non-linear programming (NLP) problem, the pseudospectral method is used with second-kind Chebyshev points.

The transformation to express the OCPs in the time interval $\tau \in [-1, 1]$ is given by
\begin{align*}
t=\frac{t_f-t_0}{2}\tau+\frac{t_f+t_0}{2}.
\end{align*}
Use Lagrange interpolation polynomials with $N+1$ points as follows: 
\begin{align*}
\mathbf X(\tau)=L_0(\tau)\mathbf x(\tau_0)+\sum^{N}_{k=1}L_k(\tau)\mathbf x(\tau_k),\\
\mathbf U(\tau)=L_0(\tau)\mathbf u(\tau_0)+\sum^{N}_{k=1}L_k(\tau)\mathbf u(\tau_k),
\end{align*}
where $\tau_0$ is the initial boundary point, and $\tau_k$, $k=1, \ldots, n$, are the $N$ collocation points, which are the zeros of the second-kind Chebyshev polynomial $U_n(\tau)$ as expressed below:
\begin{align*}
\tau_k=\cos{\frac{k}{N+1}\pi}.
\end{align*}
The weights of the Chebyshev-Gauss quadrature in this case are given by
\begin{align*}
w_k=\sqrt{1-\tau_k^2}.
\end{align*}

Then the NLP problem can be obtained as (see \cite[pp.117--118]{benson2005gauss}, \cite{fahroo2002direct})
\begin{align*}
\text{Minimize} \quad& J=\frac{t_f-t_0}{2}\sum^{N}_{k=1}w_k \| \mathbf U(\tau_k)\| \\
\text{Subject to} \quad & \frac{2}{t_f-t_0} \dot{L}_0(\tau_i)\mathbf X(\tau_0)+\frac{2}{t_f-t_0} \sum^{N}_{k=1} [\dot{L}_k(\tau_i)\mathbf X(\tau_k) - \mathbf M(\tau_i) - \mathbf U(\tau_i) ]=\mathbf 0, \\
& \| \mathbf U(\tau_k)\|-F_{\rm max} \leq 0, \\
& \mathbf S_0(\mathbf X(\tau_0))= \mathbf 0, \\
& \mathbf S_{L}(\mathbf X(\tau_f))= \mathbf 0, \\
& \mathbf S_{A}(\mathbf X(\tau_f))=\mathbf 0.
\end{align*}
where $i$ indicates the $i$th collocating point. 
The augmented cost function with the constraints combined via Lagrange multipliers is given by
\begin{align*}
\begin{split}
J_a = & \frac{t_f-t_0}{2}\sum^{N}_{k=1}w_k \| \mathbf U(\tau_k)\| - \mathbf{\Gamma}_0 \cdot \mathbf S_0(\mathbf X(\tau_0)) - \mathbf{\Gamma}_L \cdot \mathbf S_{L}(\mathbf X(\tau_f)) - \mathbf{\Gamma}_A \cdot \mathbf S_{A}(\mathbf X(\tau_f))  \\ 
& - \mu (\| \mathbf U(\tau_k)\|- F _{\rm max})  - \mathbf{\Psi} \cdot \left( \frac{2}{t_f-t_0} \dot{L}_0(\tau_i)\mathbf X(\tau_0)+\frac{2}{t_f-t_0} \sum^{N}_{k=1} [\dot{L}_k(\tau_i)\mathbf X(\tau_k) - \mathbf M(\tau_i) - \mathbf U(\tau_i)]\right).
\end{split}
\end{align*}
The remaining Karush-Kuhn-Tucker (KKT) conditions at the collocating points are given by
\begin{align*}
& \nabla_{X_k} J_a= \mathbf 0, \\
& \nabla_{U_k} J_a= \mathbf 0, \\
& \nabla_{X(\tau_0)} J_a= \mathbf 0, \\
& \mu \leq 0, \\
& \mu (\| \mathbf U(\tau_k)\|-F _{\rm max})=0.
\end{align*}
Now, the optimal control problem can be solved by using well-established NLP algorithms.

\subsection{Modified Equinoctial Elements}
The same notations as \cite{gergaud2006homotopy} are utilized here to describe the modified equinoctial orbit elements. The state variables $\mathbf X_m = (P, e_x, e_y, h_x, h_y, L) $ are defined as
\begin{align*}
P=a(1-e^2), & \quad L=\Omega+\omega+\theta, \\
e_x=e\cos{(\Omega+\omega)} , & \quad e_y=e\sin{(\Omega+\omega)},  \\
h_x=\tan \frac{i}{2} \cos\Omega, & \quad h_y=\tan \frac{i}{2} \sin\Omega,
\end{align*}
where the true longitude $L$ is the fast independent variable and the other five are slow dependent variables.

The control variables $\mathbf f = (f_r, f_t, f_h)$ in RTN (Radial-Tangential-Normal) coordinates are defined with
\begin{align*}
f_r= \mathbf F \cdot \frac{\mathbf r}{\| \mathbf r \|} , \quad
f_t=\mathbf F \cdot \left( \frac{\mathbf L}{\| \mathbf L \|} \times \frac{\mathbf r}{\| \mathbf r \|}\right) , \quad
f_h=\mathbf F \cdot \frac{\mathbf L}{\| \mathbf L \|}.
\end{align*}
Then
\begin{displaymath}
{\mathbf M}_m =\sqrt{\frac{P}{\mu}} \left( \begin{array}{ccc} 0 & 2P/W & 0 \\
\sin L & \cos L + (e_x+\cos L)/W & -Z e_y/W \\
-\cos L & \sin L + (e_y+\sin L)/W & Z e_x/W \\
0 & 0 & (\frac{C}{2} \cos L) /W \\
0 & 0 & (\frac{C}{2} \sin L) /W \\
0 & 0 & Z/W \end{array} \right),
 \mathbf N_m=\sqrt{\frac{\mu}{P}} \left( \begin{array}{c} 0 \\ 0 \\ 0 \\ 0 \\ 0 \\ \frac{W^2}{P} \end{array} \right)
\end{displaymath}
with
\begin{align*}
& W=1+e_x \cos L + e_y \sin L,\\
& Z=h_x \sin L - h_y \cos L,\\
& C=1+h^2_x+h^2_y.
\end{align*}
The system equation is given by
\begin{align}
\dot{\mathbf X}_m(t)=\mathbf M_m(t) \left( \begin{array}{lll} f_r \\ f_t \\ f_h \end{array} \right) + \mathbf N_m(t).
\label{eq:moe}
\end{align}
The minimum fuel consumption orbit transfer problem can be written as
\begin{displaymath}
(P2)  \left \{ \begin{array}{lll} \text{Minimize} &  J=\int^{t_f}_{t_0} \| \mathbf f(t) \|dt \\ 
\text{Subject to} &\dot{\mathbf X}_m(t)=\mathbf M_m(t) \mathbf f(t) + \mathbf N_m(t)\,  , \,  \forall t \in [t_0,t_f] \\
 & \| \mathbf f(t) \| \leq F_{\rm max} \, , \, \forall t \in [t_0,t_f] \\
& \mathbf X_m(t_0) \, \text{ fixed}\\
& \mathbf X_m(t_f) \, \text{ fixed}\\
\end{array}  \right.
\end{displaymath}


The optimality conditions expressed in terms of modified elements are not provided here since those presented above in Cartesian coordinates also apply to modified elements.

\section{Numerical Results}
Using an example, this section illustrates transfers from a low-Earth orbit (LEO) to a geosynchronous orbit (GEO) in terms of minimum transfer time and minimum fuel consumption. The numerical data used are from \cite[p.374]{chobotov2002orbital}, with the exception of the final time in the minimum fuel case. The initial point is given by
\begin{displaymath}
a_0=7000 \,\text{km}\, ,\, e_0=0\, ,\, i_0=28.5^{\circ}\, ,\, \Omega_0=0^{\circ} \, ,\, \omega_0=0^{\circ} \, ,\, M_0=-220^{\circ}.
\end{displaymath}
The target orbit is given by
\begin{displaymath}
a_f=42,000 \,\text{km}\, ,\, e_f=10^{-3}\, ,\, i_f=1^{\circ}\, ,\, \Omega_f=0^{\circ} \, ,\, \omega_f=0^{\circ}
\end{displaymath}
with the control constraint $F _{\rm max}=9.8\times 10^{-5}$ km/sec$^2$ and the initial time $t_0=0$. In the minimum fuel  case, the fixed final time $t_f$ is $20.17$ hr.

Canonical units are used in simulations, where 806.812 sec = 1 canonical time unit; 6378.140 km = 1 canonical distance unit; $ 9.8 \times 10^{-3}$ km/sec$^{2}$ = 1 canonical acceleration unit; and the gravitational parameter $\mu = 1$. In the following, all units are canonical unless otherwise indicated. The initial and final conditions in Cartesian coordinates are given by

\begin{displaymath}
\mathbf r_0=\left(\begin{array}{c}-0.8407\\ 0.6200\\ 0.3366\end{array}\right) , \quad \dot{\mathbf r}_0=\left(\begin{array}{c}-0.6136\\-0.6426\\-0.3489\end{array}\right),
\end{displaymath}
\begin{displaymath}
\mathbf L_T=\left( \begin{array}{c}0\\-0.0448 \\ 2.5657\end{array}  \right) , \quad \mathbf A_T=\left( \begin{array}{c}0.0010\\0\\0\end{array}  \right).
\end{displaymath}
The initial and final conditions in modified elements are given by
\begin{displaymath}
\mathbf{X}_{m0}=(1.0975 , 0 , 0 , 0.25398 , 0 , -3.8397),
\end{displaymath}
\begin{displaymath}
\mathbf{X}_{mT}=(6.5850 , 0.0010 , 0 , 0.0087 , 0 , \text{free})
\end{displaymath}
with $F _{\rm max}=0.0100$. In the minimum fuel case, the fixed final time is $t_f=90$, but in the minimum time case the final time is free, so it must first be adjusted according to the final time guess of $t_f=73$.

\begin{figure}
 \centering
\includegraphics[width=\linewidth,height=1.5in]{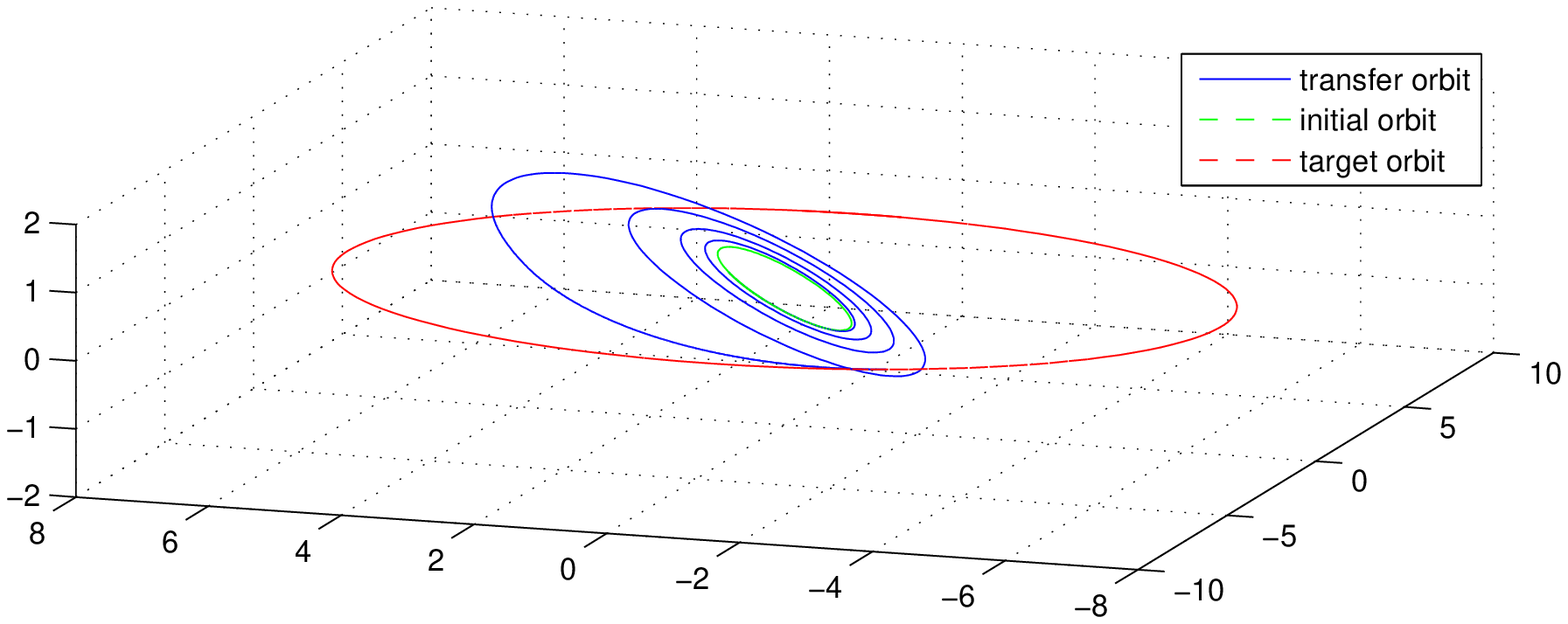}
\caption{Minimum fuel transfer (Cartesian)}
\label{fig:transfer_C}
\end{figure}

\begin{figure}
\centering
\includegraphics[width=\linewidth,height=1.5in]{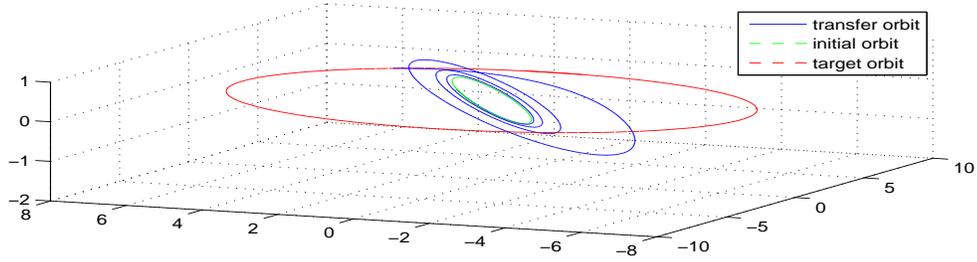}
\caption{Minimum time transfer (Cartesian)}
\label{fig:transfer_time_C}
\end{figure}

\begin{figure}
 \begin{subfigmatrix}{2}
  \subfigure[$\| \mathbf F \|$ (350 points, Cartesian)]{\includegraphics{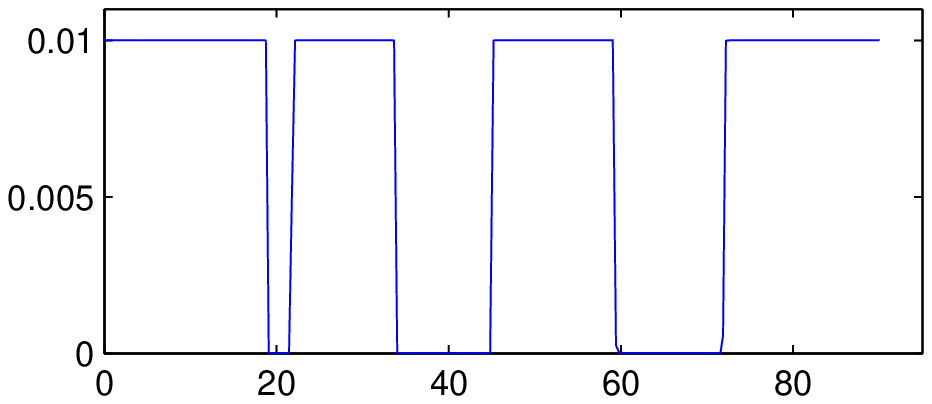}}
  \subfigure[$\| \mathbf f \|$ (300 points, Modified)]{\includegraphics{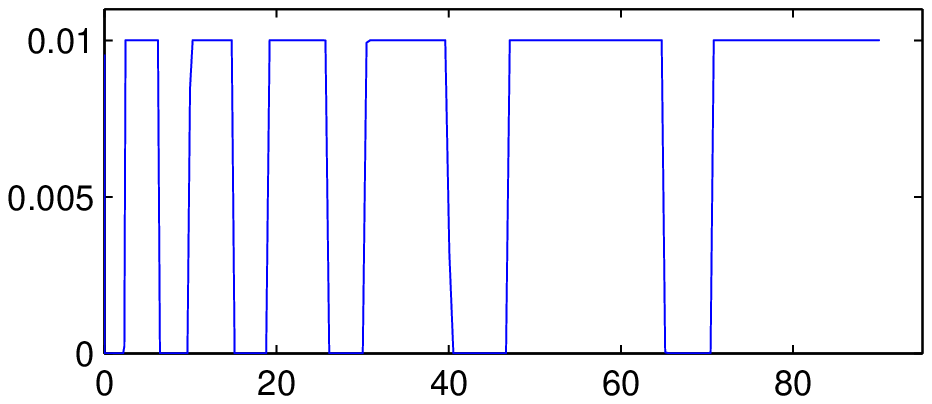}}
 \end{subfigmatrix}
 \caption{Comparison of control force magnitude time history of min fuel}
 \label{fig:force comparison C}
\end{figure}

\begin{figure}
 \begin{subfigmatrix}{2}
  \subfigure[$\| \mathbf L \|$]{\includegraphics{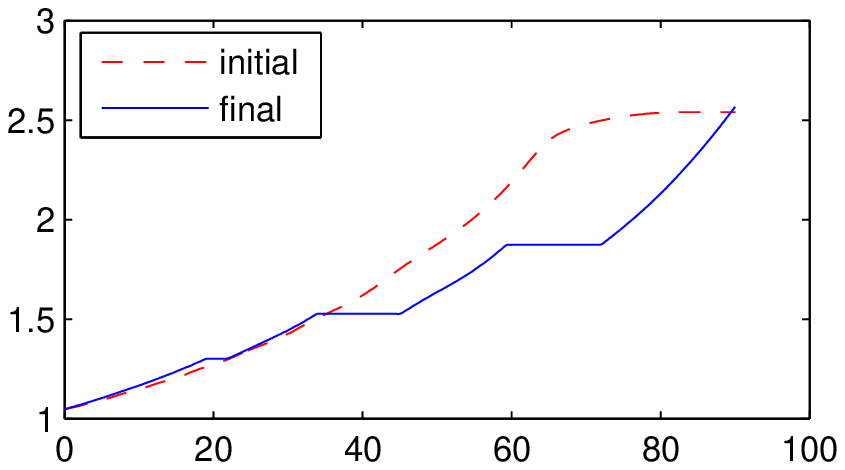}}
  \subfigure[$\| \mathbf A \|$]{\includegraphics{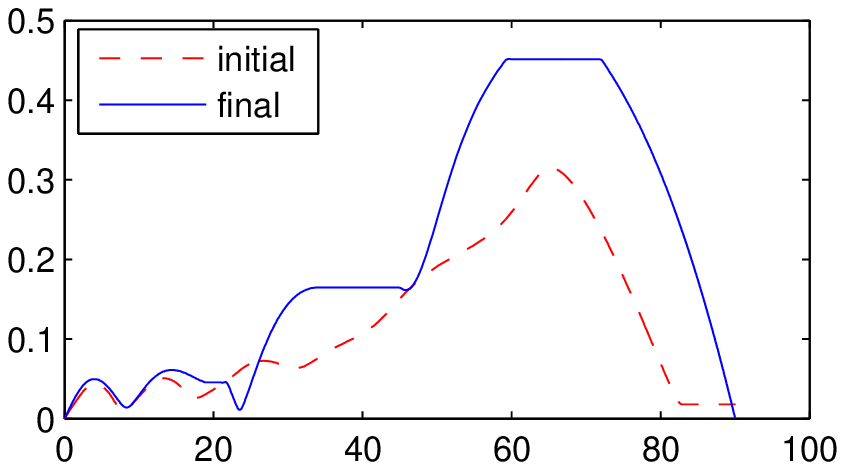}}
 \end{subfigmatrix}
 \caption{$\| \mathbf L\|, \| \mathbf A\|$ time history of min fuel (Cartesian)}
 \label{fig:LA C}
\end{figure}

\begin{figure}
 \begin{subfigmatrix}{3}
  \subfigure[$r_1$]{\includegraphics{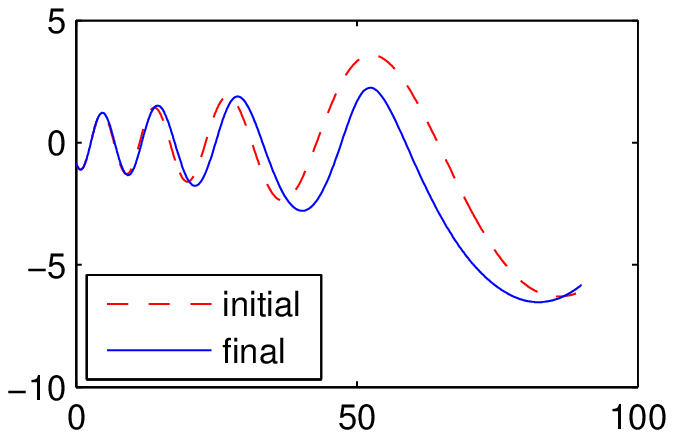}}
  \subfigure[$r_2$]{\includegraphics{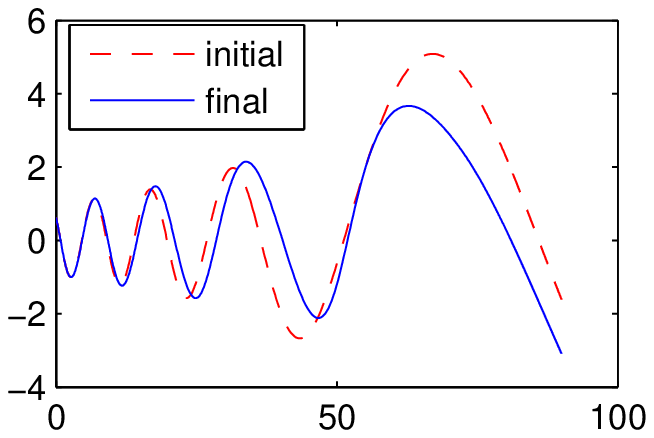}}
  \subfigure[$r_3$]{\includegraphics{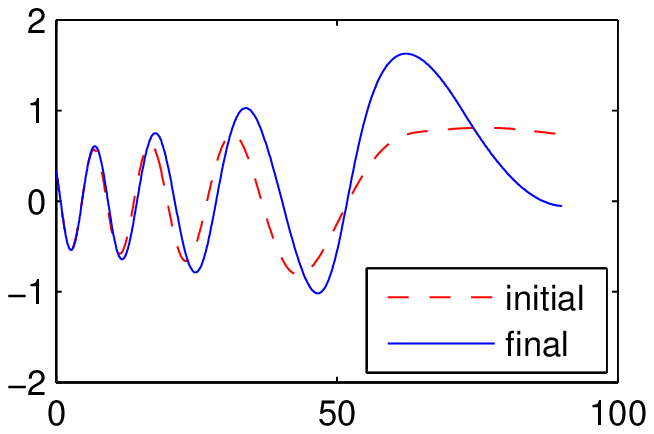}}
  \subfigure[$\dot r_1$]{\includegraphics{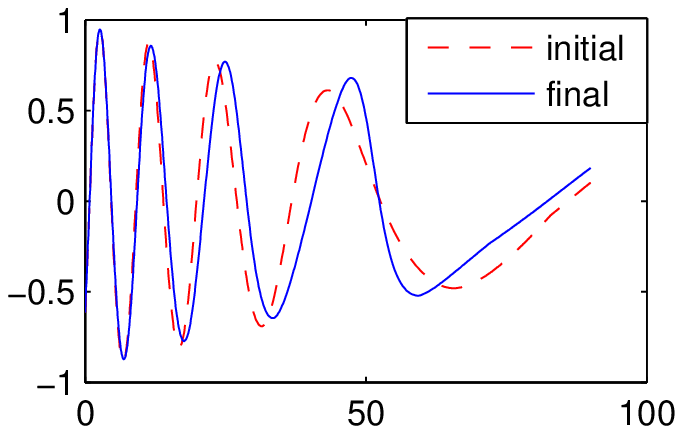}}
  \subfigure[$\dot r_2$]{\includegraphics{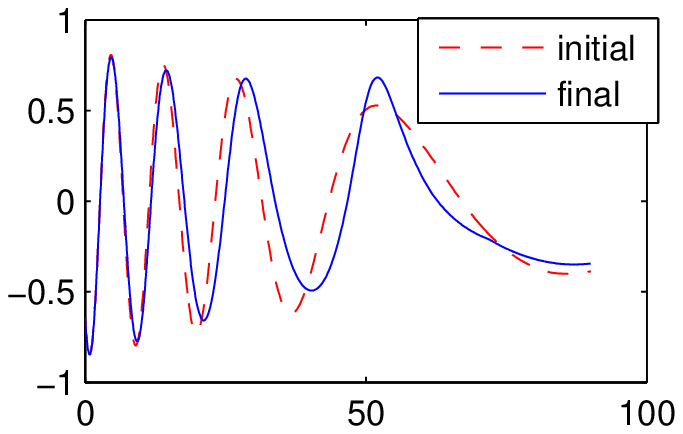}}
  \subfigure[$\dot r_3$]{\includegraphics{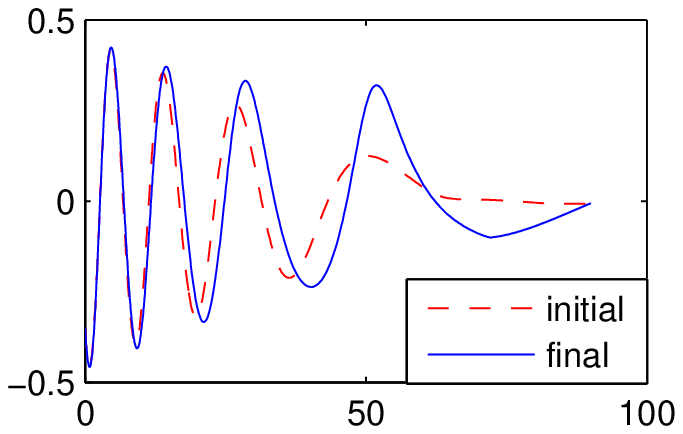}}
 \end{subfigmatrix}
 \caption{$(\mathbf r, \dot{\mathbf r})$ time history of min fuel (Cartesian)}
 \label{fig:r_dr_C}
\end{figure}

\begin{figure}
 \begin{subfigmatrix}{3}
  \subfigure[$F_1$]{\includegraphics{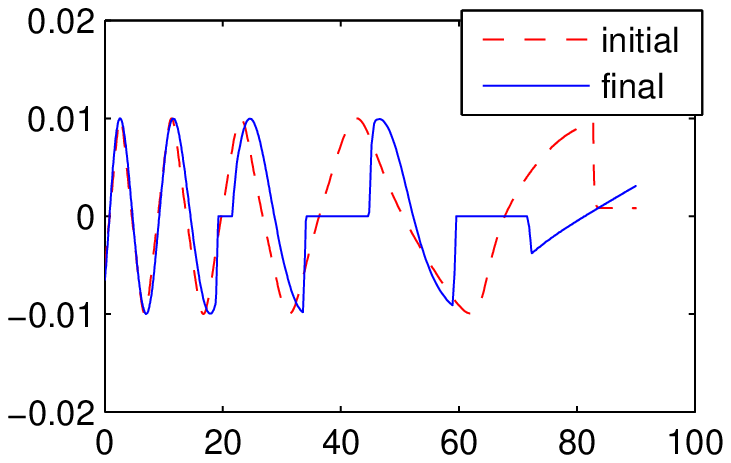}}
  \subfigure[$F_2$]{\includegraphics{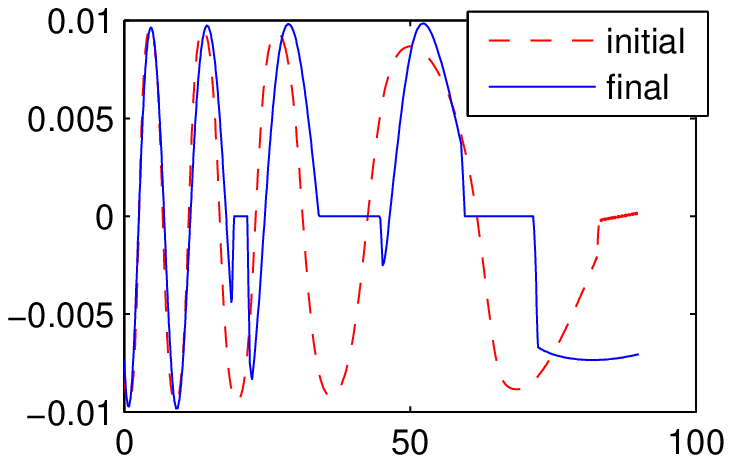}}
  \subfigure[$F_3$]{\includegraphics{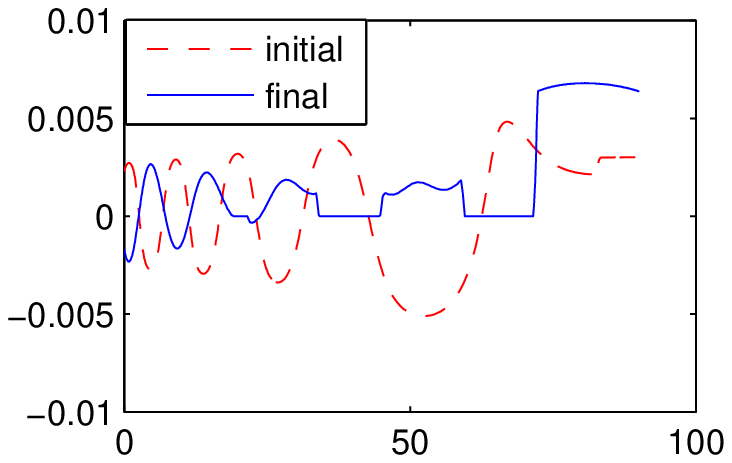}}
 \end{subfigmatrix}
\caption{Control force time history of min fuel (Cartesian)}
\label{fig:force_C}
\end{figure}

The CCM controller given in (\ref{eq:control}) is applied with the weighting $k=2$ and the function $f=1$ to obtain a transfer trajectory that provides an initial guess for optimization. Then TOMLAB/PROPT is utilized together with the pseudospectral method and SNOPT solver to optimize the trajectory on the MATLAB platform. To compare the different coordinate systems, all the collocating points are located in one phase. The optimal results are listed in Table \ref{Table:C}, and the minimum fuel and minimum time transfer trajectories are shown in Fig. \ref{fig:transfer_C} and Fig. \ref{fig:transfer_time_C}, respectively. For the minimum fuel case, Fig. \ref{fig:force comparison C} compares the force magnitude between the two coordinate systems presented above.  Fig. \ref{fig:LA C} and Fig. \ref{fig:r_dr_C} show the time history of $\| \mathbf L \|, \| \mathbf A \|$  $\mathbf r$,  and $\dot{\mathbf r}$, and Fig. \ref{fig:force_C} details the time history of the control force in the Cartesian coordinates system. In the time history plots, the dotted lines represent the initial guesses and the solid lines represent final trajectories.

\begin{table}
\centering
    \begin{tabular}{  c  c  c  c  c  c c }
    \hline
     & & Transfer    & Switch  & $\Delta V$ & Collocating  &Optimality \\ [2pt] 
     & Basis &  Time (hr)   &  Times &  (km/sec) & Points & Condition \\ [2pt] \hline
  \multicolumn{7}{c}{Method in this paper} \\  [2pt]   \hline
     & Cartesian &  20.1703  & 6 & 4.9911 & 300 & Satisfied \\  [2pt] 
     &  &  20.1703  & 6 & 4.9938 & 350 & Satisfied \\ [2pt]   \cline{2-7}
    Min Fuel &  & 20.1703 & 11 & 4.9945 & 300 & Not Satisfied \\  [2pt] 
     & Modified & 20.1703 & 12 & 5.1788 & 350 & Not Satisfied \\  [2pt] 
     &  & 20.1703 & 12 & 5.1772 & 400 & Not Satisfied \\ [2pt]  \hline
     Min Time & Cartesian & 16.1362 & 0 & 5.6873 & 300 & Satisfied \\ [2pt]  \cline{2-7}
      & Modified  & 18.2205 & 0 & 6.4212 & 300 & Satisfied \\ [2pt]   \hline
      \multicolumn{7}{c}{Chobotov's Result in \cite[pp.374-376]{chobotov2002orbital}} \\ [2pt]   \hline
     Min Time & Equinoctial  & 16.2845 & 0 & 5.7451 & * & * \\ [2pt] \hline
    \end{tabular}
    \caption{\label{Table:C}Comparisons of optimal results in different basis} 
 \end{table}

\begin{figure}
\includegraphics[width=\linewidth,height=1.5in]{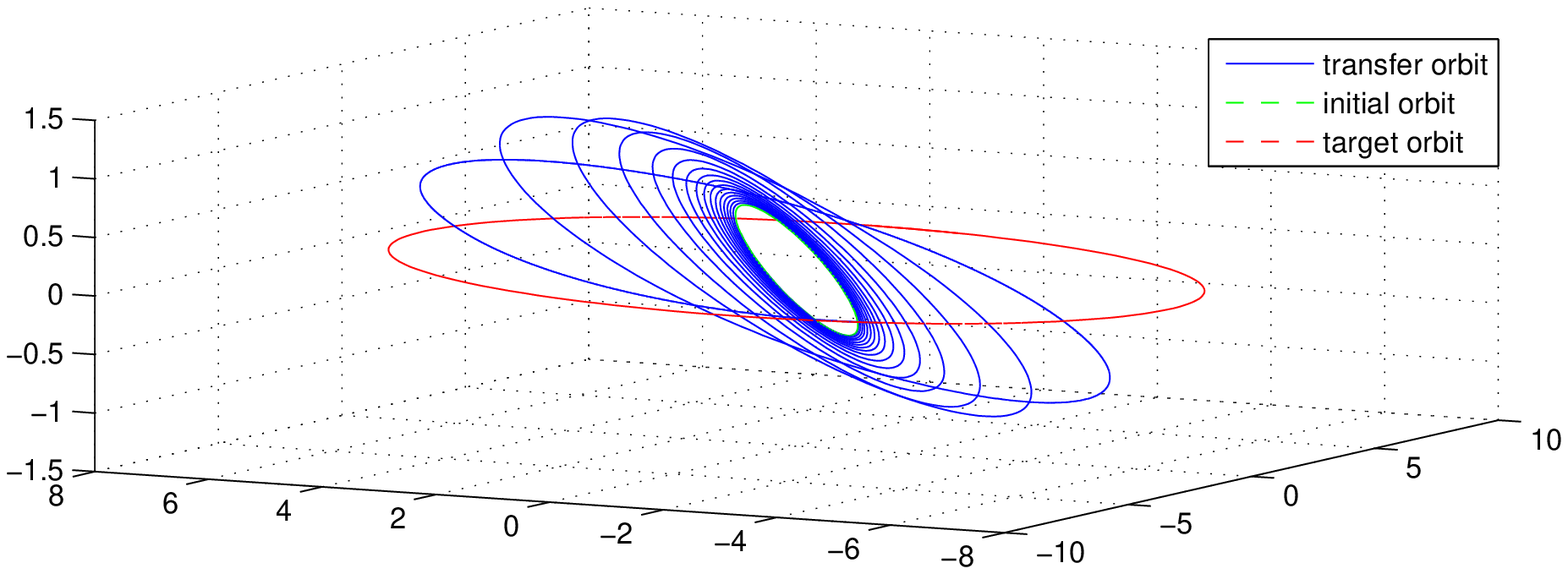}
\caption{Minimum fuel transfer (Cartesian, $ F_{\rm max} /4$ )}
\label{fig:transfer_C_over4}
\end{figure}

\begin{figure}
 \begin{subfigmatrix}{3}
  \subfigure[$r_1$]{\includegraphics{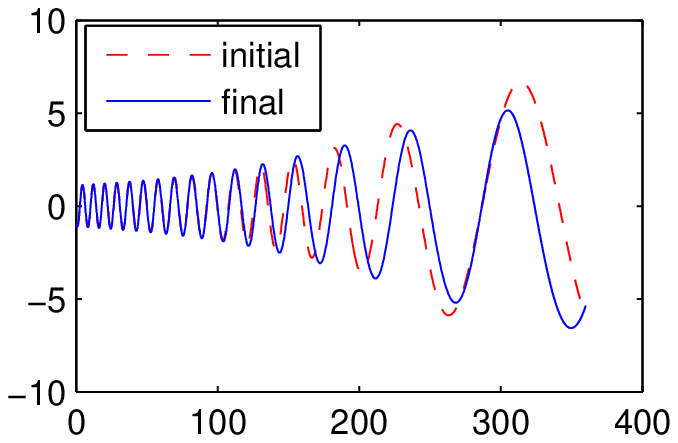}}
  \subfigure[$r_2$]{\includegraphics{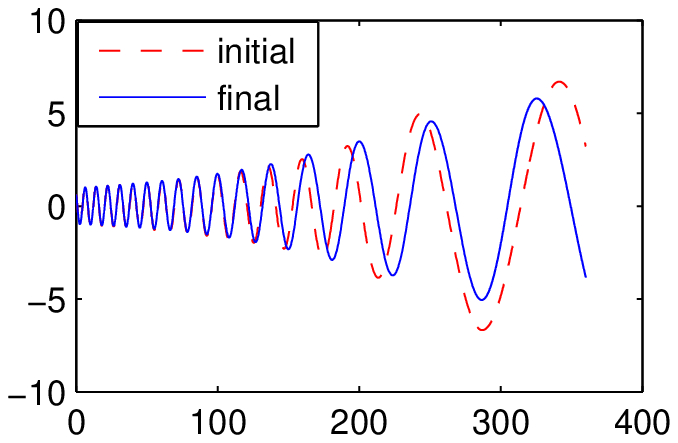}}
  \subfigure[$r_3$]{\includegraphics{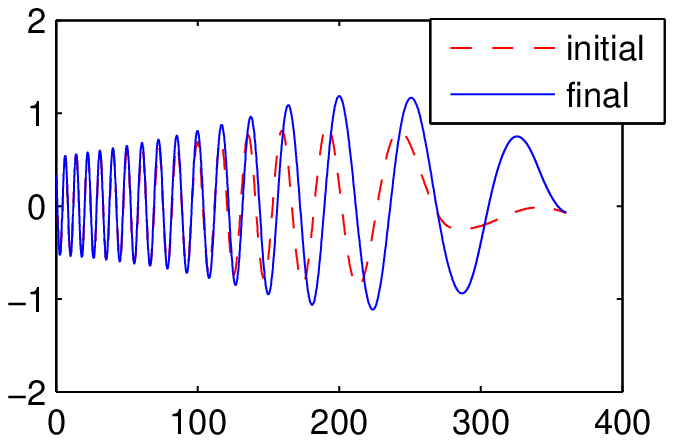}}
  \subfigure[$\dot r_1$]{\includegraphics{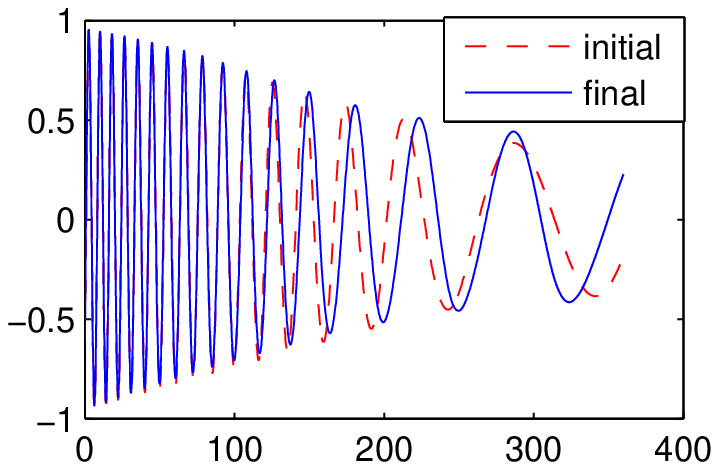}}
  \subfigure[$\dot r_2$]{\includegraphics{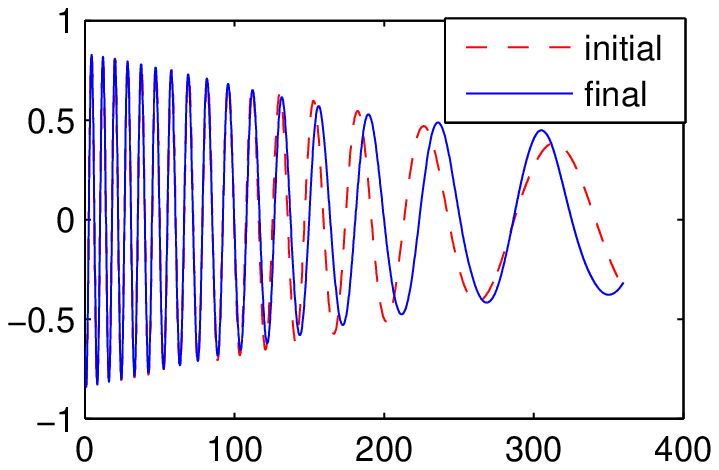}}
  \subfigure[$\dot r_3$]{\includegraphics{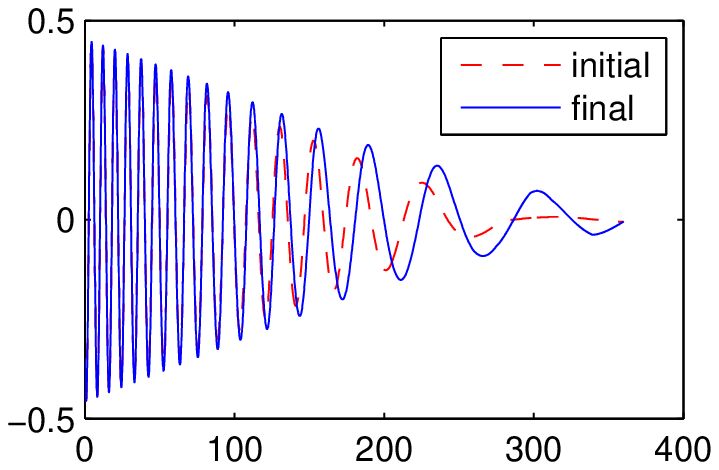}}
 \end{subfigmatrix}
 \caption{$(\mathbf r, \dot{\mathbf r})$ time history of min fuel (Cartesian, $ F_{\rm max} /4$ )}
 \label{fig:r_dr_C_over4}
\end{figure}

Narrower control constraints are then set by $\frac{1}{2} F _{\rm max}$ and $\frac{1}{4} F_{\rm max}$ to simulate the optimal fuel trajectory with the fixed final time $2t_f$ and $4t_f$. The optimal results are listed in Table \ref{Table:C_lower}. The final transfer trajectory and the time histories of the state variables with $\frac{1}{4} F_{\rm max}$ are shown in Fig. \ref{fig:transfer_C_over4} and Fig. \ref{fig:r_dr_C_over4}.

\begin{table}
\centering
    \begin{tabular}{  c  c  c  c  c  c  c c  }
    \hline
     &  & Control  & Transfer    & Switch  & $\Delta V$  & Collocating  & Optimality  \\ [2pt] 
     & Basis &  Constraints &  Time (hr)   &  Times &  (km/sec) &  Points &  Conditions \\ [2pt] \hline
     &  &  $F _{\rm max} $ & 20.1703 & 6 & 4.9938 & 350 &  \\ [2pt] 
     Min Fuel & Cartesian &  $ F_{\rm max}/2$ & 40.3406  & 11 & 5.1185 & 400 & Satisfied \\  [2pt] 
       &  & $F_{\rm max} /4$ &  80.6812  & 31 & 5.1747 & 400 & \\ [2pt]  \hline
    \end{tabular}
\caption{\label{Table:C_lower}Optimal results of low-thrust minimum fuel problem}

 \end{table}
 
In the minimum time case, the optimal result of transfer time obtained using the Cartesian coordinates system is similar to Chobotov's in \cite{chobotov2002orbital}, while the optimized result obtained using modified elements is far from optimal. The mathematical description of the modified elements is an approximation from the equinoctial elements found by the stroboscopic method. Therefore, it fails to approximate the real optimal solution in the presented example. However, the equinoctial elements are not a good choice of a coordinate system for direct optimization since Kepler's equation remains to be solved.

In the minimum fuel consumption case, the results from the use of Cartesian coordinates and those from the use of the modified elements with 300 collocating points in Table \ref{Table:C} are similar to each other. Even still, all the optimality conditions  are not satisfied in the case where the modified elements are used. Fig. \ref{fig:force comparison C} shows six switching times of $\| \mathbf F \|$ in the Cartesian coordinates system with 350 collocating points and eleven switching times of $\| \mathbf f \|$ in the modified elements with 300 collocating points. This implies that using modified elements requires more collocating points and higher order polynomials to approximate the trajectory. Obtaining the optimal solution in this case increases the difficulty of computation noticeably. Even when modified elements are applied to find an optimal solution close to that which was obtained via Cartesian coordinates, this method would again fail to approximate the real optimal solution since the optimized control force sequences would be quite different from those in Cartesian coordinates. Moreover, problematic results can be identified with 350 and 400 collocating points from modified elements as illustrated in Table \ref{Table:C}. However, by using  Cartesian coordinates that contain weaker nonlinearity than orbital elements, obtaining the real optimal solution through polynomial approximation can be achieved with ease.

The time histories of $(\mathbf r, \dot{\mathbf r})$ in Fig. \ref{fig:r_dr_C} and Fig. \ref{fig:r_dr_C_over4} indicate the importance of the CCM transfer as the initial guess of trajectory optimization in the Cartesian coordinates system. While the NLP solver is quite robust, without a supported Lyapnov-based initial guess, it would still be difficult for the solver to move those collocating points onto the optimal trajectory, due to the periodically changing sign and rapidly changing value of the state variables in the Cartesian coordinates system \cite[p.52]{conway2010spacecraft}.

Besides the strong nonlinearity in the equations of motion, the use of modified elements has several numerical disadvantages because its system, which is the right-hand side of equation (\ref{eq:moe}), is only an approximation. For a given low-thrust transfer trajectory, the time histories calculated from the right-hand side of the differential equation (\ref{eq:moe}) can not match the real time histories of the variables' derivatives on left-hand side of (\ref{eq:moe}). In particular, for stiff problems -- as with long duration cases -- where the trajectories are lacking high accuracy from the ODE solver, the error between the two sides of the equation will be large.

\section{Conclusions}

In this paper, a simple and effective approach has been proposed by employing the pseudospectral method,  nonlinear programming, and the Chang-Chichka-Marsden transfer controller. Solutions to the free-injection minimum fuel  consumption  and  minimum time transfer problems have successfully  been obtained by using Cartesian coordinates supported by a Chang-Chichka-Marsden transfer trajectory.  The Chang-Chichka-Marsden transfer, as an initial guess, has revealed its usefulness for overcoming the numerical difficulties associated with the direct optimization caused by strong oscillation of state variables in the Cartesian coordinates system. The use of orbital elements increases the difficulty of optimization and fails to provide the optimal solution by the nonlinear programming and pseudospectral methods. However, by utilizing Cartesian coordinates instead of orbital elements, the optimal solutions are easy to obtain.

Two main advantages arise when Cartesian coordinates are used in the direct trajectory optimization. First, this option is the simplest way to describe the two-body system accurately without the optimal solution-related problems posed by singularity or approximation. Second, the weaker nonlinearity makes it easier to obtain the optimal solution, mainly due to the numerical approximation method used.

Future research on longer duration transfer described in Cartesian coordinates can be performed by developing a specific numerical method for quick and accurate computation with fewer collocating points. Additionally, the oblateness of the Earth and the shadow effect shall be taken into account.

\section{Acknowledgments}
This research was carried out in the Department of Applied Mathematics, University of Waterloo, and the School of Astronautics, Harbin Institute of Technology. The authors are indebted to Per Rutquist in TOMLAB for the advice  on the use of the  optimization software PROPT.

\newpage

\begin{tabular}{ p{1.6in} p{3.2in}  }

\vspace{0pt}
\includegraphics[width=0.9\linewidth]{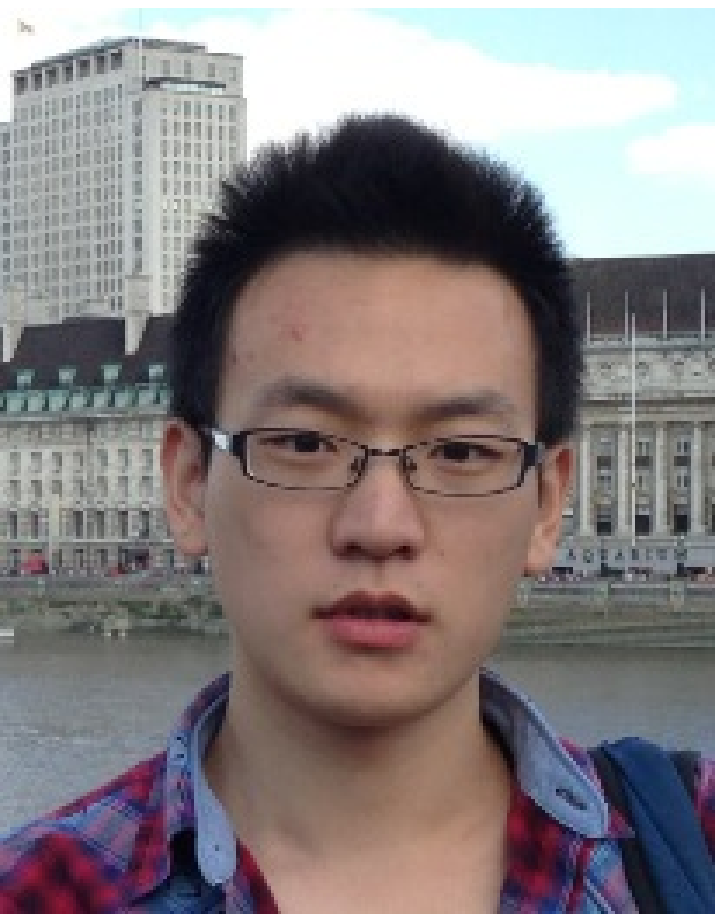} 
&
\vspace{0pt}
\textbf{Hantian Zhang} is a current undergraduate student at Department of Aerospace Engineering and Mechanics, School of Astronautics, Harbin Institute of Technology. He expected to receive his B.Eng.  from Harbin Institute of Technology in 2014. He worked as a research assistant with Prof. Dong Eui Chang at Department of Applied Mathematics in University of Waterloo from January to March in 2013. His research interests include orbital dynamics and control, trajectory optimization, nonsmooth mechanics, and dynamical systems.
\\
\vspace{20pt}
\includegraphics[width=0.9\linewidth]{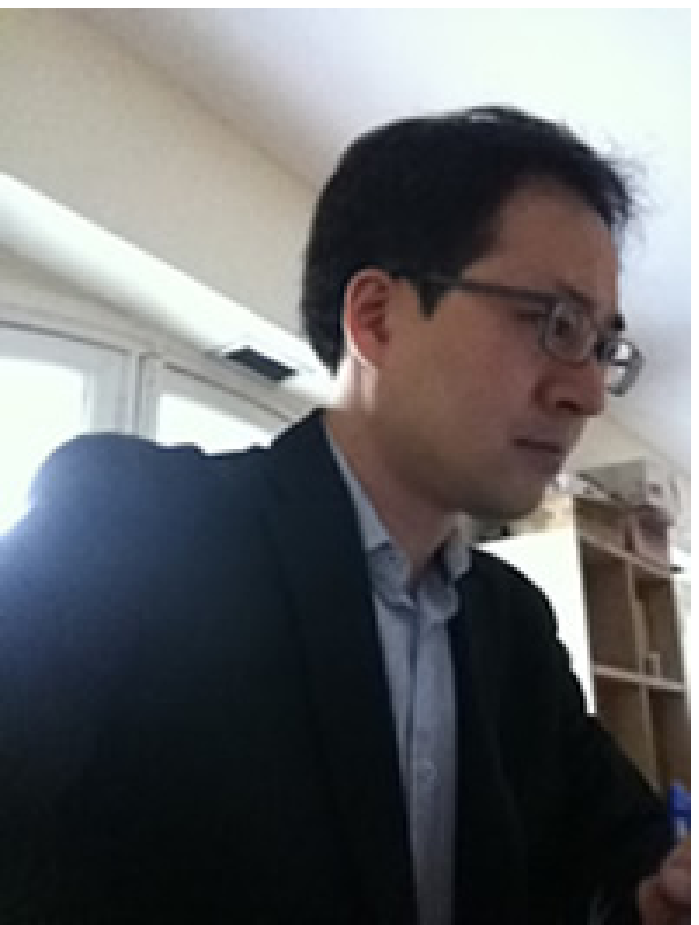}
&
\vspace{20pt}
\textbf{Dong Eui Chang} received his B.S. and M.S. from Seoul National University in 1994 and 1997, and his Ph.D. in Control \& Dynamical Systems (CDS) from the California Institute of Technology (Caltech). He is currently an Associate Professor in the  Department of Applied Mathematics at the University of Waterloo. His research interests lie in control, mechanics and various engineering applications.
\\
\vspace{20pt}
\includegraphics[width=0.9\linewidth]{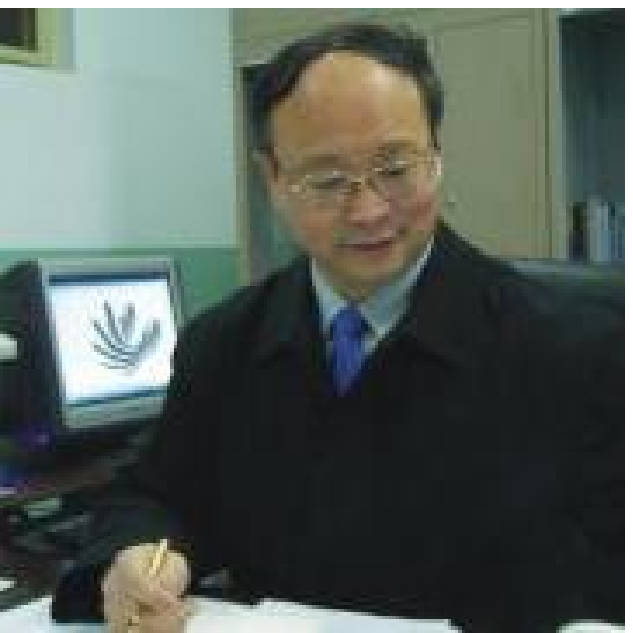}
&
\vspace{20pt}
\textbf{Qingjie Cao} received his B.S. and M.S. in Mathematics from Qufu Normal University in 1982 and 1985, and his Ph.D. in Mechanical Engineering from Tianjin University in 1993. Before he joined Harbin Institute of Technology, he served as a Professor at School of Mathematics in Shandong University, and Research Fellows in University of Aberdeen and University of Liverpool. Currently, he is a Professor at School of Astronautics in Harbin Institute of Technology, and the director of the Centre for Nonlinear Dynamics Research. He introduced the SD oscillator and SD attractor in mechanical systems.
\end{tabular}

\end{document}